\documentclass[letterpaper, 10 pt, conference]{IEEEtran}  

\IEEEoverridecommandlockouts                              

\usepackage{url}
\usepackage{cite}
\usepackage{hyperref}
\usepackage[utf8]{inputenc}
\usepackage[english]{babel}
\let\labelindent\relax
\usepackage{enumitem}
\usepackage{graphics} 
\usepackage{epsfig} 
\usepackage{mathptmx} 
\usepackage{times} 
\usepackage{amsmath} 
\usepackage{amssymb}  
\usepackage{mathtools}
\usepackage{lipsum}
\usepackage{cancel}
\newtheorem{theorem}{Theorem}
\newtheorem{definition}{Definition}
\newtheorem{observation}{Observation}
\newtheorem{corollary}{Corollary}

\newtheorem{remark}{Remark}

\title{\LARGE \bf
Design of Incentive Mechanisms Using Prospect Theory to Promote Better Sell-back Behavior among Prosumers
}

\author{Diptangshu Sen$^a$ and Arnob Ghosh$^b$
\thanks{$^{a}$ Diptangshu Sen is an undergraduate alumnus of the Department of Mechanical Engineering, Indian Institute of Technology Delhi, India.
        {\tt\small dsen.research@gmail.com}}%
\thanks{$^{b}$ Arnob Ghosh is an Assistant Professor at the Department of Mechanical Engineering, Indian Institute of Technology Delhi, India. He is currently also a  research associate at the Control and Power Engineering Group of the Electrical and Electronic Engineering Department at the Imperial College, London. 
        {\tt\small arnob.ghosh@mech.iitd.ac.in, arnob.ghosh@imperial.ac.uk}}%
}

\begin{document}

\maketitle
\thispagestyle{empty}
\pagestyle{empty}

\begin{abstract}
Users can now give back energies to the grid using distributed resources. Proper incentive mechanisms are required for such users, also known as prosumers, in order to maximize the sell-back amount while maintaining the retailer's profit.   However, all the existing literature considers expected utility theory (EUT) where they assume that prosumers maximize their expected payoff in the face of uncertainty. We consider prospect theory (PT) which models the {\em true} behavior of humans under uncertainty. We show that in a day-ahead contract pricing mechanism, the actual optimal value of contract and the sell-back amount may be smaller compared to the one computed by EUT. We also propose a lottery-based mechanism and show that such a mechanism can increase the sell-back amount while increasing the retailer's savings compared to day-ahead contract pricing.  
\end{abstract}

\section{Introduction}
\subsection{Motivation}
The distribution side of the traditional power grid is rapidly changing. Consumers now have access to renewable energy sources like solar cells and windmills. Hence, consumers can now produce or give back energy to the grid. The rise of such consumers (also denoted as 'prosumers') are critical for realizing the goal of reducing carbon emissions. Further, in the advent of excess demand at a certain location, those energies can be used to mitigate the supply and demand imbalance.  However, without proper incentivization mechanism the prosumers may not give back energies to the grid.

One of the most popular ways of incentivization is through payment at a retail-rate for every unit of energy contributed by the prosumer. In this {\em net-metering} approach, the prosumers are paid for their {\em net} contribution at the retail rate. Though it is simple, it leads to a ``utility death spiral" for retailer (\cite{death}). For example, there is a loss in transmission, so the rate at which the prosumers are compensated is higher than the retail rate. Hence, in order to maintain the profit, the retailer needs to increase the retail rate which in turn will increase the compensation rate for prosumers and by extension, the retail price. Thus, the consumers will end-up paying a high price.  Further, the retailer traditionally procures energy at the wholesale price which is smaller than the retail price. Some policymakers propose real-time pricing for feed-in energy, however, real-time prices are often volatile. If the wholesale price is very low, the prosumers would not have any incentive to give back. Further, high volatile prices render the prosumer's decision volatile which makes the retailer's decision to procure energies volatile. In order to reduce the volatility, a day-ahead contract pricing approach has been considered. In this price mechanism, a day-ahead price is announced to the prosumers, the prosumers then commit to a contract. The prosumers are compensated at this day-ahead price for their committed amounts. In real-time, if a prosumer fails to fulfill his commitment, he would incur a penalty for the shortfall amount. 

Though such a price mechanism is attractive, one key component of developing such a price mechanism is anticipating how a prosumer would respond to the uncertainty inherent in renewable energy generation. Should the prosumer commit a large amount in order to get a guaranteed day-ahead return while risking the shortfall amount in the real-time? Or, should the prosumer commit a smaller amount in order to avoid any penalty in the real-time? Researchers have answered these questions using expected utility theory (EUT) where they consider that prosumers will maximize the expected utility. However, prosumers are human beings, and it is well established that they are not risk-neutral. Thus, expected utility theory may not give the correct estimate of the feed-in energy from the prosumers. 

\subsection{Goals}
We seek to answer the following questions:
\begin{itemize}
\item What is the optimal amount of energy given back by the prosumers for a given day-ahead contract price? The answer to such a question is useful to the retailer, since the retailer needs to plan based on the amount fed back by the prosumers. Further, it will settle the debate whether EUT can provide accurate predictions for the total sell-back energy. 

\item Does there exist a mechanism which can incentivize the users to give back more while maintaining the profit of the retailer? If we can identify a simple mechanism which can achieve both objectives, it will be a {\em win-win} scenario for both the retailer and the prosumers. 
\end{itemize}
\subsection{Our Approach}
In order to answer the first question, we resort to{\em prospect utility theory} which models human decision-making behavior under uncertainty. In day-ahead contract pricing, the prosumer makes a decision in two phases: firstly on a day-ahead basis, he chooses the optimal contract amount that maximizes his payment and in real-time, he chooses what amount of energy to feed back such that his overall utility is maximized. Prosumer's utility has two components, one of them is the utility/convenience derived from consumption and while the other is the utility derived from the payment he receives. The first-stage decision differs in prospect utility theory compared to expected utility theory, since the consumer needs to consider the uncertainty of real-time realization of the renewable energies. We show that the committed contract  is smaller compared to the committed contract predicted by EUT when the penalty price is larger than a threshold. Our analysis also shows that the real-time sell-back amount is smaller than the amount computed using EUT because of the risk-averse behavior of prosumers. 

In order to answer the second question, we investigate a fixed lottery-based incentive scheme for sell-back where the chance of a prosumer winning the lottery is proportional to the sell-back energy contributed by him. Although the probability of winning is very small, prospect theory postulates that human beings have a tendency of overestimating small probabilities. This along with the fact that there is no penalty,  makes the lottery scheme very lucrative to the prosumers . Empirically, we show that the lottery scheme achieves significant improvements over contract pricing in terms of sell-back amounts. Our analysis also reveals that the retailer saves more in the lottery-based scheme compared to the day-ahead contract based scheme.

\subsection{Related Literature}
In recent times, a lot of research effort has gone into developing energy sharing and management systems for smart grids involving 'prosumers'(\cite{cerquides2015designing,skopik2012novel,zafar2018prosumer,wang2015review}). Game-theoretic techniques are very popular in the analysis of distributed energy systems(\cite{lee2015distributed}, \cite{saad2011coalitional}, \cite{saad2011noncooperative}, \cite{alskaif2018gamification}, \cite{chen2019buy}).  Besides feeding energy back to the grid, other kinds of efficient energy sharing techniques have also been studied in literature(\cite{neagu2020new,saad2011noncooperative,saad2011coalitional,ghosh2020penalty,menniti2013demand}). 

In \cite{chiu2014optimal}, a time-dependent day-ahead pricing model for pricing sell-back energies is proposed. Others(\cite{bitar2010optimal}, \cite{bitar2011role} and \cite{bitar2019marginal}) investigate optimal day-ahead contract pricing for incentivising WPPs(Wind power producers) who sell back. Some papers also consider how contracts would be affected in presence of storage facilities(\cite{bitar2011role} and \cite{bitar2019marginal}). However, all these papers consider EUT to model uncertainty while computing the optimal contract amount. In contrast, we use prospect theory which closely emulates the behavior of humans. 

Some papers consider prospect theory for modeling energy sell-back behavior of prosumers under day-ahead contract pricing as well as random future pricing(\cite{el2016prospect}, \cite{el2017managing}). \cite{el2016prospect} uses a non-cooperative game to find the optimal day-ahead contract structure for all prosumers and shows that prospect theory contracts tend to be smaller. However, we find closed-form expressions for prosumer's optimal contract by accounting for his real-time decision. We find conditions under which PT contracts are smaller than EUT contracts. We also propose a lottery scheme and show that it can produce higher sell-backs than contract pricing. 

Lottery mechanisms for demand response problems in smart grids have been studied by some papers(\cite{schwartz2014demand} and \cite{kim2019optimal}), while coupon-based systems for demand response have also been explored(\cite{ming2020prediction}). However, the above papers do not consider prospect theory while designing the mechanism. In contrast, we provide a theoretical understanding of the lottery mechanism using prospect theory, and compute the optimal sell-back amount for a given lottery scheme. Further, we evaluate the impact of the lottery scheme on net sell-back amount and retailer's savings. 

\section{Background: Prospect Theory}
EUT postulates that when presented with multiple risky choices, a rational human always chooses the option which offers the highest expected payoff. However, it has been shown experimentally in \cite{kahneman2013prospect} that humans react differently to losses and gains. Further, humans tend to overestimate events which occur with lower probabilities compared to the events which occur with higher probabilities. Prospect theory tries to fix some of the inconsistencies in EUT(\cite{kahneman2013prospect}). In the following, we describe the prospect theory in detail. 
\subsection{Value Function}
We define $v(y)$ as the perceived value function associated with a gain by amount $y$ from a reference point. The reference point is often taken as a guaranteed deterministic return or $0$.
\begin{itemize}
    \item Human beings are usually loss-averse. Therefore, for $y > 0$, $-v(-y) > v(y)$. This means that a person is more disappointed by a loss of amount $y$ than she is happy about a gain of amount $y$ (\cite{kahneman2013prospect}).
    \item It can be proved that an exponential form of $v(y)$ is most suitable and is considered throughout literature. For $y > 0$, $v(y) = y^{\eta}$, and  $v(-y) = \lambda(-y)^{\beta}$ (proposed by \cite{kahneman2013prospect}, proved by \cite{al2008note}).
    \item From the loss-averse property, we must have $-\lambda(-y)^{\beta} > y^{\eta}$. Plugging $y = 1$, we have $\lambda > 1$. $\lambda$ has been empirically found to be close to 2 (\cite{al2008note}).
    \item Using the property of diminishing sensitivity with increase in value of $|y|$ (\cite{tversky1992advances}), it can be shown that $\eta, \beta \leq 1$. 
\end{itemize}
It has been shown recently that for persons with homogenous preferences, $\eta = \beta$ (\cite{al2008note}). If losses and gains are not too significant with respect to the reference value, $\eta = \beta = 1$ is a good approximation. Throughout this paper, we assume that $\eta=\beta=1$.

\subsection{Decision Weight}
Let $\pi(\cdot)$ denote the function that maps probability to decision weights. We denote $\pi(\cdot)$ as the weight function. Note that $\pi(\cdot)$ is not a probability measure. $\pi(\cdot)$ is monotonically increasing in probability $q$ and we have $\pi(0) = 0$ and $\pi(1) = 1$. Since humans tend to overestimate lower probability events, $\pi(q) > q$ when $q$ is close to $0$. Further, humans underestimate higher probability events, $\pi(q)<q$ when $q$ is close to $1$. One of the widely used decision weight functions is given by the following
\begin{equation}\label{decisionweight}
    \pi(q) = \exp(-(-\log q)^{\gamma})
\end{equation}
where $0 < \gamma < 1$. This form was first proposed by \cite{prelec1998probability}. We use $\gamma = 0.5$ throughout this paper. Note that $\pi(q)$ is concave initially and then convex with an inflexion point at $q = \frac{1}{e}$.

\subsection{Perceived value of Prospects}
From \cite{kahneman2013prospect}, a prospect $(y_1, q_1; y_2, q_2;...y_n,q_n)$ is a contract that yields outcome $y_i$ with probability $q_i$ where $q_1+q_2+...q_n=1$. For simplification, we can omit null outcomes and use the simplified notation $(y, q)$ to represent the prospect. Let $\textbf{V}(y, q)$ be the prospect theory utility associated with the said contract. Then
\begin{equation*}
    \textbf{V}(y, q) =\sum_{i} v(y_i)\pi(q_i) 
\end{equation*}
If all outcomes in a prospect have a common deterministic component, we can remove it from the prospect in the \textit{editing phase}.  Let us consider a prospect of the form $(y_1 + c, q; y_2 + c, 1-q)$ where $c$ is the fixed gain or loss. Then we have
\begin{equation*}
    \begin{multlined}
    \textbf{V}(y_1 + c, q; y_2 + c, 1-q) = c + \textbf{V}(y_1, q; y_2, 1-q) \\= c + v(y_1)\pi(q) + v(y_2)\pi(1-q)
    \end{multlined}
\end{equation*}

\subsection{Cumulative Prospect Theory}
In 1992, a modification was proposed to the original theory which came to be known as \textit{Cumulative Prospect Theory (CPT)}(\cite{tversky1992advances}). The primary contribution of CPT is that it preserves first order stochastic dominance and can be extended to cases with infinitely many outcomes. From CPT, the utility associated with the random variable $Y$ and its distribution function $F$ as 
\begin{equation}\label{cpt}
    \begin{multlined}
    U(y,F) = \int_{-\infty}^{0} v(y)\frac{d}{dy}(\pi(F(y)))dy \\ + \int_{0}^{\infty}v(y)\frac{d}{dy}(-\pi(1-F(y)))dy
    \end{multlined}
\end{equation}
Here, $\pi(\cdot)$ is the weight function corresponding to the distribution function $F$. Note that while facing loss, humans become risk-seeking, hence, they tend to think what is the probability that the loss can be less than or equal to $y$. While facing gain, humans tend to think what is the probability that the gain can be greater than or equal to $y$, hence, they become risk-averse. The CPT formalizes the utility according to the above psychology.

\section{System Model}
We consider a model where there is a community of $N+n$ consumers and an utility company/retailer $R$ who supplies electricity to the community. Out of $N+n$ consumers, $n$ are prosumers. The time is slotted into $T$ periods that span 24 hours. Each consumer $i$ is characterized by a parameter $\omega_k^{(i)}$ which denotes his \textit{willingness for demand}\footnote{From the consumer's payoff function in Eq. 4, we can deduce that the maximum demand consumer $i$ is willing to consume is given by $\frac{\omega_i^{(k)}}{\alpha}$ which is realized when price $p_k = 0$. Even at zero price, he will not consume infinite amounts as it does not increase his convenience any further. Since the maximum theoretical demand depends directly on $\omega_i^{(k)}$, we call it the ''willingness for demand'' parameter of consumer $i$.} in the $k^{th}$ period. All consumers derive some comfort from the consumption of energy which is not always tangible. We monetize that comfort by defining the \textit{convenience function} $C(x_k^{(i)}, \omega_k^{(i)})$ where $x_k^{(i)}$ is the amount of energy consumed by the $i^{th}$ consumer in the $k^{th}$ period. We have
\begin{equation}\label{eq:utility}
    C(x_k^{(i)}, \omega_k^{(i)}) = 
    \begin{cases} 
      \omega_k^{(i)} x_k^{(i)} - \alpha\frac{(x_k^{(i)})^{2}}{2} & x_k^{(i)} \leq \frac{\omega_k^{(i)}}{\alpha} \\
      \frac{(\omega_k^{(i)})^2}{2\alpha} &  x_k^{(i)} \geq \frac{\omega_k^{(i)}}{\alpha}
    \end{cases}
\end{equation}
The use of convenience functions is well documented in literature (\cite{samadi2010optimal}, \cite{fahrioglu1999designing} and \cite{fahrioglu2001using}), sometimes with varying nomenclature. $\alpha$ is a  constant. Observe that the convenience function has several nice properties, apart from being continuous and twice differentiable. It has a fixed point at the origin because zero consumption must lead to zero convenience. The first derivative is positive indicating that higher consumption leads to higher convenience. It is also concave meaning that marginal convenience derived per unit consumption decreases as consumption increases. And finally, the convenience function saturates beyond a threshold. Quadratic convenience functions are quite common in the smart grid and economic theory literature.

The consumer $i$ also pays a price $p_k$ for consuming energy from the grid at time period $k$. The consumer $i$'s payoff is
\begin{align}\label{eq:payment}
    C(x_k^{(i)}, \omega_k^{(i)})-p_kx_k^{(i)}
\end{align}

For a prosumer $j$, let $s_k^{(j)}$ be the renewable energy generation in the $k^{th}$ period. We assume that in some period $k$, prosumers produce renewable energy in excess and do not need to purchase from the grid. However, it is to be noted that in other periods, prosumers might not produce enough renewable energy and have to purchase from the grid. Let $z_k^{(j)}$ denote the amount fed back to the grid by the $j^{th}$ prosumer in said period. Clearly, $0 \leq z_k^{(j)} \leq s_k^{(j)}$. Therefore, $s_k^{(j)}-z_k^{(j)}$ is the actual amount consumed by him in that period. Hence, the convenience function can be modified to the following form: 
\begin{equation}\label{conv}
    \begin{multlined}
    C(z_k^{(j)}, \omega_k^{(j)}) = \\
    \begin{cases} 
      \omega_k^{(j)} (s_k^{(j)}-z_k^{(j)}) - \alpha\frac{(s_k^{(j)}-z_k^{(j)})^{2}}{2} & (s_k^{(j)}-z_k^{(j)}) \leq \frac{\omega_k^{(j)}}{\alpha} \\
      \frac{(\omega_k^{(j)})^2}{2\alpha} &  (s_k^{(j)}-z_k^{(j)}) > \frac{\omega_k^{(j)}}{\alpha}
    \end{cases}
    \end{multlined}
\end{equation}
Equation (\ref{conv}) is central to the rest of our analysis. Please note that henceforth, we use a simplified notation where we drop the prosumer index from the superscript for ease of expression. Note that we have not considered the temporal correlation among the utility functions across multiple time periods.  

\section{Day-Ahead Contract Pricing}
In this scheme, the retailer announces the sell-back price of the $k^{th}$ period, $p_s$ on a day-ahead basis. The prosumers are paid based on the committed amount, however, if they fail to satisfy the committed amount they have to pay a penalty price $p_e$ for every unit of short-fall. Note that the existing literature so far only considers the setting where prosumers are distributed generators for the day-ahead price mechanism and do not consider their own convenience while committing sell-back amounts. We consider a more generalized setting where a prosumer takes his own convenience into account while choosing sell-back contracts.

Let us consider a prosumer $P$ with a willingness for demand parameter $\omega_k$ in period $k$. Suppose, $P$ signs a contract in the day-ahead market, promising to sell-back an amount $C_k$ in the $k^{th}$ period. However, in real time, he has a renewable energy generation $s_k$ and he decides to sell-back amount $z_k$ which may not be equal to $C_k$. Thus, the prosumer's decision is in two time-scales, in the day-ahead he commits an amount to sell-back, later when the renewable energies are realized, he decides the actual amount to sell-back based on the penalty. Thus, in the day-ahead time, the prosumer must also consider the decision in real-time for computing an optimal amount to commit. First, we describe the prosumer's decision in real-time. Subsequently, we describe the prosumer's decision in the day-ahead market using CPT. We also characterize the difference in committed amounts between CPT and EUT.

\subsection{Prosumer's decision in real-time}
In real-time, $P$ has to choose the sell-back amount $z_k$. Note that the choice of $C_k$ has already been made on a day-ahead basis and hence it is a constant when choosing $z_k$. Note that in the real-time, everything is realized, thus, there is no uncertainty, hence, we do not need to consider CPT. The prosumer's utility is real-time is:
\begin{equation}
    U(z_k) = C(z_k, \omega_k) + p_s C_k - p_e(C_k-z_k)^{+}
\end{equation}
$C_k$ is already given, the prosumer decides $z_k$. $(C_k-z_k)^{+}$ denotes the shortfall amount on which penalty is accrued.
 Since $U(z_k)$ is not differentiable at $z_k = C_k$, we proceed with a case-by-case analysis. Note that when $z_k > C_k$, $z_k^{*} = s_k - \frac{\omega_k}{\alpha}$. When $z_k < C_k$, $z_k^{*} = s_k - \frac{(\omega_k-p_e)}{\alpha}$. Thus, we obtain--

\begin{observation}\label{realtime}
Let, $z_1 = s_k - \frac{\omega_k}{\alpha}$ and $z_2 = s_k - \frac{(\omega_k-p_e)}{\alpha}$. The optimal decision in the real-time:
\begin{itemize}
    \item If $C_k < z_1$, $z_k^{*} = z_1$
    \item Else if, $C_k > z_2$, $z_k^{*} = z_2$
    \item Else, $z_k^{*} = C_k$
\end{itemize}
\end{observation}
Note that when the generation is low, the prosumer reduces his own consumption in order to avoid penalties. Also when $p_e$ is high, in order to avoid penalty, he sells a higher amount by reducing his own consumption. 

\subsection{Prosumer's decision on a day-ahead}
On a day-ahead basis, the prosumer needs to choose the optimal contract amount $C_k$.  Note that the penalty price should always be greater than the sell-back rate. Otherwise, a prosumer can always choose a very large contract, without selling back anything, yet he can attain a high payoff. When choosing the contract, the aim of the prosumer is to maximize the payment according to prospect theory because only that component of his utility depends on the contract.We denote this payment by $E_p(C_k)$.

 In order to use (\ref{cpt}), we need to define gains/losses for each choice of sell-back amount $z_k$. Let $y_k = z_k-C_k$. Note that after editing of prospects (removal of fixed component $p_s C_k$), we are only left with probabilistic loss terms.
 We define two new parameters:
 \begin{definition}
  $z_{2,min}$ is the minimum value assumed by $z_2$ and it occurs when renewable energy generation in real-time $s_k$ attains its minimum value $s_{min}$. 
 
  $z_{2,max}$ is the maximum value assumed by $z_2$ which is attained when $s_k = s_{max}$. 
 \end{definition}
 $C_k$ cannot exceed $z_{2,max}$ because from Observation 1, the highest sell-back is given by $z_{2,max}$. Hence, if $C_k$ is greater than $z_{2,max}$, a prosumer will certainly incur a loss of $p_e(C_k-z_{2,max})$. Therefore, under CPT the valuation of prosumer $P$ is : 
\begin{equation*}
    \begin{multlined}
    E_p(C_k) = p_s C_k + p_e \int_{z_{2,min}}^{C_k}v(z_k-C_k)\frac{d}{dz_k}(\pi(F_z(z_k)))dz_k
    \end{multlined}
\end{equation*}
where $F_z(\cdot)$ is the CDF of $z_k$.
Setting $\frac{dE_p(C_k)}{dC_k} = 0$ and changing variables from $z_k$ to $s_k$, we have the following result:
\begin{theorem}\label{thm1}
The optimal contract amount predicted by prospect theory is given by $C_k^{*} = F_s^{-1}(\pi^{-1}(\frac{p_s}{\lambda p_e})) + \frac{(p_e - \omega_k)}{\alpha}$.\\
\end{theorem}
The proof for Theorem \ref{thm1} has been relegated to the online preprint version (\cite{sen2020design}).
The corresponding EUT expression for $C_k^{*}$ can be easily found by setting $\lambda = 1$ and observing that $\pi^{-1}(x) = x$. This gives us $C_k^{*} = F_s^{-1}(\frac{p_s}{p_e}) + \frac{(p_e-\omega_k)}{\alpha}$ for an EUT prosumer. Note that as penalty price $p_e$ increases, the expression in the first term decreases, however, the expression in the second term increases. Hence, the variation of the contract amount with the penalty price is not clear. Intuitively, the contract amount should decrease as the penalty price increases. However, note from Observation 1 that when the penalty price is high, the sell-back amount increases in order to avoid hefty penalty, which in turn makes the prosumer commit to a higher amount. If we have a distributed energy generator instead of a prosumer, then $\omega_k = 0$ and the contract amount decreases with an increase in penalty price, as expected.  

The next result provides the condition under which the actual contract amount is smaller compared to the value computed by EUT.\\
\begin{corollary}\label{thm:cpt}
$p_s < \frac{p_e}{e} \implies F_s^{-1}(\pi^{-1}(\frac{p_s}{\lambda p_e})) < F_s^{-1}(\frac{p_s}{p_e} )$. 
\end{corollary}
This result indicates that if the penalty price exceeds a threshold, the sell-back amount committed by a prosumer is actually smaller compared to the one predicted by EUT. The intuition behind this result is that prosumers are risk-averse when faced with uncertainty and end up taking more conservative choices compared to those predicted by EUT. However, the result also indicates that if the penalty price is not very large, the prosumer's committed amount in the day-ahead market may also be higher than what is computed by EUT. 

\begin{corollary}\label{thm:sell_back}
$C_{EUT} \geq C_{PT}$ $\implies$ $z_{EUT}^{*} \geq z_{PT}^{*}$.  \end{corollary}
This result indicates that if EUT predicts a higher contract amount compared to CPT, then the actual sell-back amount as computed in CPT  will be smaller than the one computed in EUT. The above result suggests that if a retailer computes prosumer sell-back amounts using EUT, there may be significant errors in the computation. Combining the results from Corollaries \ref{thm:cpt} and \ref{thm:sell_back}, it is easy to discern that when the penalty price is high, the sell-back computed in EUT is smaller than the actual amount computed in CPT.\\

\begin{remark}
Some papers also suggest a penalty price for selling back in excess of the amount committed in the day-ahead. Our analysis can easily be extended to the above case. 
\end{remark}

\section{Fixed Lottery Mechanism}
In this section,we propose an incentive scheme. In the next section, we, empirically, show that such an incentivization scheme  increases the sell-back amount and the savings of the retailer. 

In this mechanism, there are no contracts or penalties. Instead, the retailer announces a lottery of fixed prize $R$. A prosumer can win the entire amount. The chance of a prosumer winning the lottery is directly proportional to the energy contribution made by him. This ensures that the mechanism is fair. At the end of the period, the retailer decides the winner of the lottery. We investigate the actual sell-back amount in this mechanism, and whether this mechanism mandates higher sell-back while resulting in a higher saving for the retailer. 

 If prosumer $P$ contributes an amount $z_k$ in the $k^{th}$ period, his chance of winning the lottery is given by $q = m z_k$ where $m$ is a constant. $m$ needs to be chosen carefully such that $\sum_{j=1}^{n} q_j \leq 1$ where $j \in J$, the index set for prosumers. (1-$\sum_{j}q_j$) is the probability that the lottery is won by none and the lottery amount might be carried forward to the next round. The utility for said prosumer is given by:
\begin{equation}\label{lottery_eqn}
  \begin{multlined}
  U(z_k) = \omega_k(s_k-z_k)-\frac{\alpha}{2}(s_k-z_k)^{2} + \textbf{V}(R, q) \\
  = \omega_k(s_k-z_k)-\frac{\alpha}{2}(s_k-z_k)^{2} + R\pi(m z_k)
  \end{multlined}
\end{equation}
This follows because $\textbf{V}(R, q) = v(R)\pi(mz_k) = R\pi(mz_k)$. The optimal sell-back amount $z_k^{*}$ is given by solving:
\begin{equation}\label{eq:firstorder}
    \begin{multlined}
    U^{'}(z_k) = -\omega_k + \alpha(s_k-z_k) + Rm\pi^{'}(m z_k) = 0
    \end{multlined}
\end{equation}
We use the form of $\pi(\cdot)$ defined in (\ref{decisionweight}). When $p$ is small, $\pi(p)$ is concave and we have a unique maximum for (\ref{lottery_eqn}). The optimal sell-back amount $z_k^{*}$ is obtained by solving the equation:
\begin{equation*}
    -\omega_k + \alpha(s_k-z_k) + R\gamma\frac{\exp(-(-\log mz_k)^{\gamma})(-\log(mz_k))^{\gamma-1}}{z_k} = 0
\end{equation*}
There is a unique solution to the above equation denoted as $z$, hence, the optimal value will also be unique. \\
\begin{theorem}
A higher lottery always earns a higher sell-back, i.e., $R_2 > R_1$ $\implies$ $z_{2,i}^* > z_{1,i}^*$ $\forall$  $i = 1,\ldots,n$ where $z_{2,i}^*$ ($z_{1,i}^*$, resp.) is the solution of (\ref{eq:firstorder}) for $R=R_2$ ($R_1$, resp.).\\
\end{theorem}
Hence, the above result indicates that as the lottery prize increases, the sell-back amount from a prosumer also increases. The detailed proof can be found in \cite{sen2020design}.
\section{Numerical Experiments}\label{results}
\subsection{Setup}
We consider that there are $N = 7,500$ consumers in the community. For consumer $i$, $\omega_i$ is drawn uniformly from the interval $[3, 7]$ independently from other consumers. From (\ref{eq:payment}), it can be easily discerned that in response to the retail price $p_k$, consumer $i$ consumes $\frac{\omega_i-p_k}{\alpha}$.  The total consumption of the consumers is denoted as the \textit{base demand $B$}. Throughout our analysis, we use a retail rate $p_k$ at $Rs. 1.5$/unit.

In addition to the 7,500 consumers, we have $n = 2,500$ prosumers who have access to in-house renewable energy and they do not need to buy energy from the grid in the given period. They commit to a day-ahead market at the rate $p_s$. We assume that prosumers are slightly high-end customers. Thus, for prosumer $i$, $\omega_i$ is uniformly drawn from an interval  $[4, 7]$.  For each prosumer $j$, there is uncertainty around the renewable energy generation $s_j$. We consider that $s_j > \frac{\omega_j}{\alpha}$ for all $j = 1,\ldots,n$. Precisely, $s_j = \frac{\omega_j}{\alpha} + r_j$ where $r_j$ is the uncertain component. $r_j$ is assumed to be drawn uniformly from the interval $[0, 0.5]$ for all $j$. Since no prosumer is able to sell-back more than 10 units theoretically, we set $m = \frac{0.1}{n}$ for the lottery model. This ensures that $\sum_{j}p_j \leq 1$ always. 

Let $Z$ denote the total sell-back by the prosumers. Given $B$ is the base demand, $D = B-Z$ is the net demand of the grid (amount that retailer needs to purchase from the wholesale market). Since the retailer needs to purchase lower amounts of energy from the wholesale marketplace, sell-back generates savings. The savings can be calculated as follows :
\begin{equation*}
    Savings = \mathbb{Q}(B) - \mathbb{Q}(D) - I
\end{equation*}
where $I$ is the cost incurred by the retailer for the payment made to the prosumers in various incentivisation schemes (lottery or contract-pricing), $\mathbb{Q}(\cdot)$ represents the quadratic cost to serve the demand. Costs are always convex and of the form $\mathbb{Q}(x) = ax + bx^2$. $a$ and $b$ are constants. For all our numerical analysis, we use $a = 1$ and $b = 2 \times 10^{-5}$. 

In case of lottery incentivization, $I = R$ where $R$ is the fixed lottery amount. For day-ahead contract pricing, the retailer has to make payments to prosumers for the amounts they sell-back at the rate $p_s$. Therefore, $I = \sum_{j=1}^{n} p_s C_j - p_e (C_j - z_j)^{+}$ where $p_e$ denotes the penalty per unit shortage from the stipulated contract. 

\subsection{Results}
\subsubsection{EUT vs PT contracts}
\begin{figure}[!ht]
     \centering
     \includegraphics[width = 0.48\textwidth]{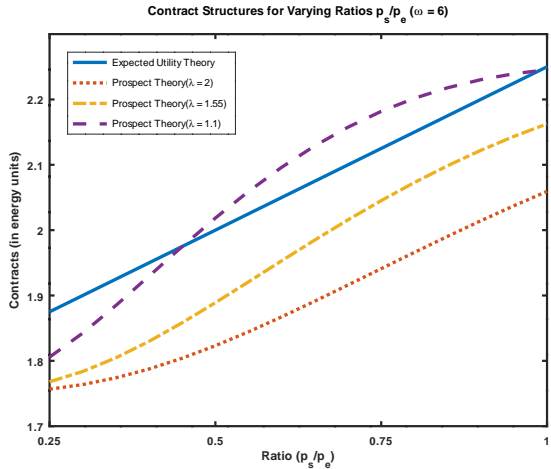}
     \vspace{-0.3in}
     \caption{Comparison between EUT and PT contracts}
     \label{EUTvsPT_comparison}
     \vspace{-0.15in}
\end{figure}

We investigate if there is a difference in contract amounts computed by EUT and CPT. When the penalty rate is fixed, contract amounts are expected to increase with increase in the sell-back rate. This trend is consistent across EUT and PT results. However, due to the risk-aversive nature of the prosumers, PT contracts are usually lower than EUT contracts (as observed in Fig. \ref{EUTvsPT_comparison}). Therefore, as $\lambda$ (parameter indicating prosumer risk-aversiveness) increases, contract amounts sharply decrease. For the most popular choice of $\lambda=2$, the contracts computed by CPT are at least $10\%$ smaller than those computed by EUT.

\subsubsection{Sell-back \& Retailer Savings}
\begin{figure}[!ht]
     \centering
     \includegraphics[width = 0.48\textwidth]{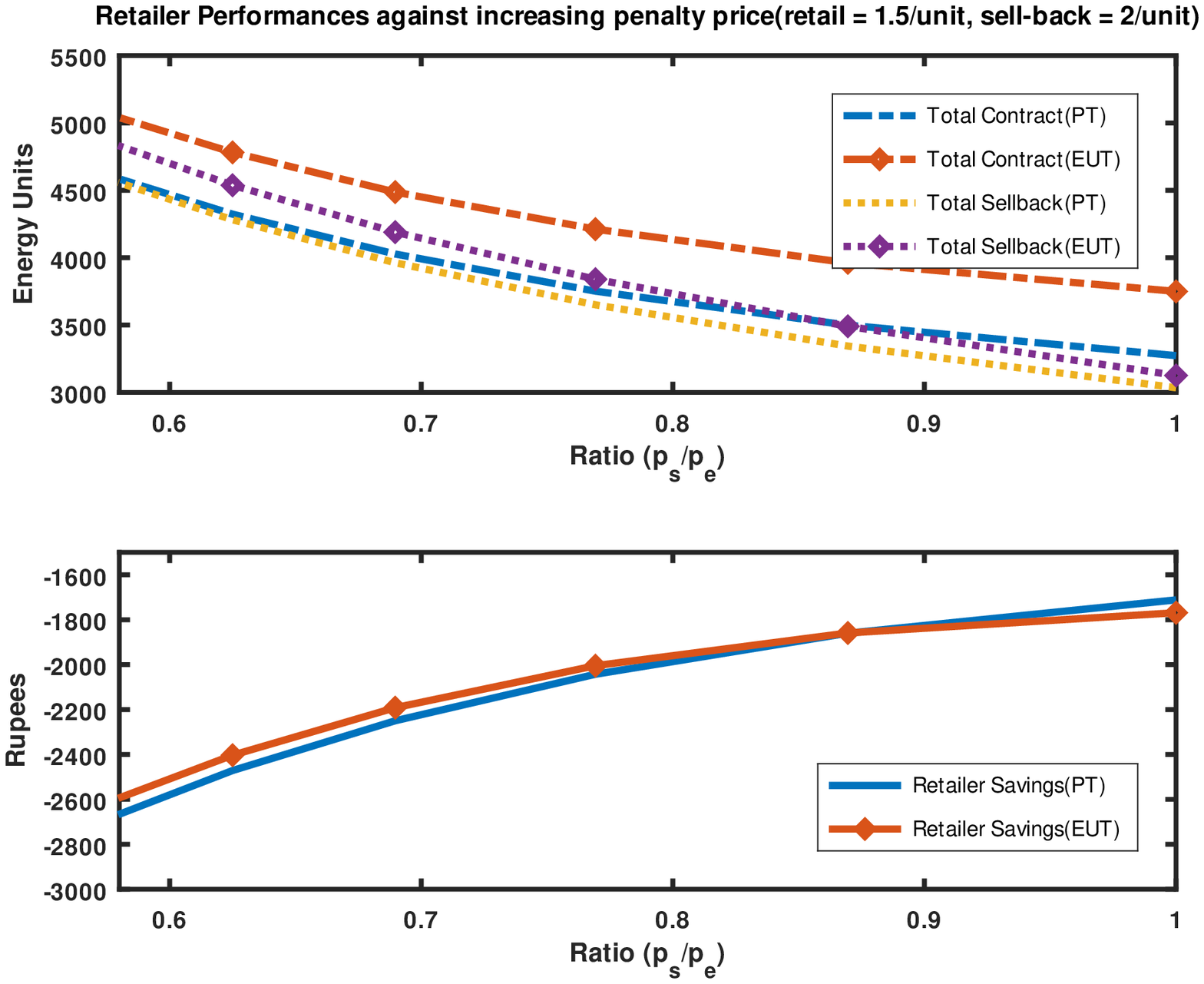}
     \vspace{-0.3in}
     \caption{Variation of metrics against varying penalty price (sellback price $>$ retail price)}
     \label{penalty_vary}
     \vspace{-0.15in}
\end{figure}

\begin{figure}[!ht]
     \centering
     \includegraphics[width = 0.48\textwidth]{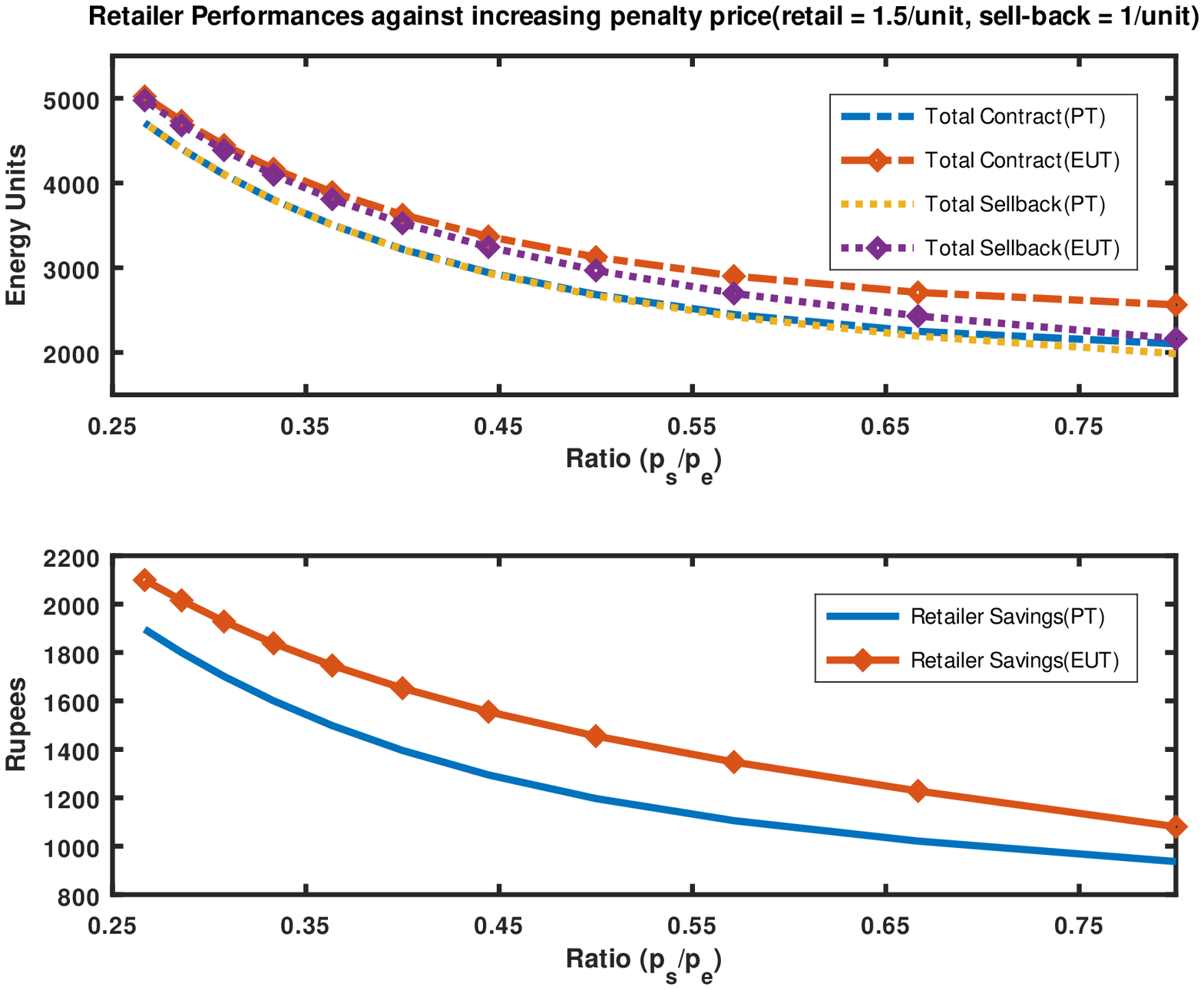}
     \vspace{-0.3in}
     \caption{Variation of metrics against varying penalty price (sellback price $<$ retail price)}
     \label{penalty_vary2}
     \vspace{-0.15in}
\end{figure}

\begin{figure}[!ht]
     \centering
     \includegraphics[width = 0.48\textwidth]{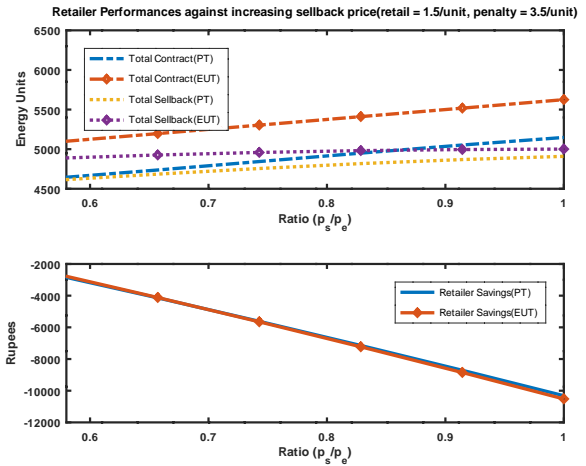}
     \vspace{-0.3in}
     \caption{Variation of metrics against varying sell-back price}
     \label{sellback}
     \vspace{-0.15in}
\end{figure}

It is observed that the actual sell-back realized in real-time is always less than the total signed contract amount in the day-ahead market (Fig. \ref{penalty_vary}, \ref{penalty_vary2} and \ref{sellback}). This is expected because prosumers have no incentive to sell-back more than what is stipulated in their contracts. Sell-backs predicted by PT are also smaller than their EUT counterparts similar to Fig.~\ref{EUTvsPT_comparison}. 

In the first set of experiments (Fig. \ref{penalty_vary}), the sell-back contract rate $p_s$ is kept fixed at $Rs. 2$/unit (higher than retail price) while the penalty price $p_e$ is varied from $Rs. 2$/unit to $Rs. 3.5$/unit. As the ratio $\frac{p_s}{p_e}$ approaches 1 (i.e., penalty price decreases), contracts and sell-back amounts gradually decrease as suggested by Theorem~\ref{thm1}. However, the retailer's saving is negative in this scenario.

When we fix $p_s$ at $Rs.1$/unit (lower than the retail price), the trends for contracts and sell-back amounts are similar to Fig.\ref{penalty_vary}. Observe that a lower sell-back rate than retail rate (Fig.\ref{penalty_vary2}) leads to positive retailer savings compared to negative savings in Fig.\ref{penalty_vary}. Thus, this price mechanism may be more beneficial to the prosumers. However,  a sell-back rate lower than the retail rate also makes sell-back a less lucrative option for prosumers which leads to lower sell-back amounts compared to Fig.\ref{penalty_vary}. 

Subsequently, we fix the penalty price $p_e$ at $Rs. 3.5$/unit and increase the sell-back price gradually (see Fig. \ref{sellback}). As expected, a higher sell-back rate incentivises prosumers to choose higher contracts and also sell-back more. However, the retailer incurs larger costs because she has to make larger payments to prosumers and therefore her savings decrease sharply. 

\subsubsection{Performance of a fixed lottery scheme}
\begin{figure}[!ht]
     \centering
     \includegraphics[width = 0.48\textwidth]{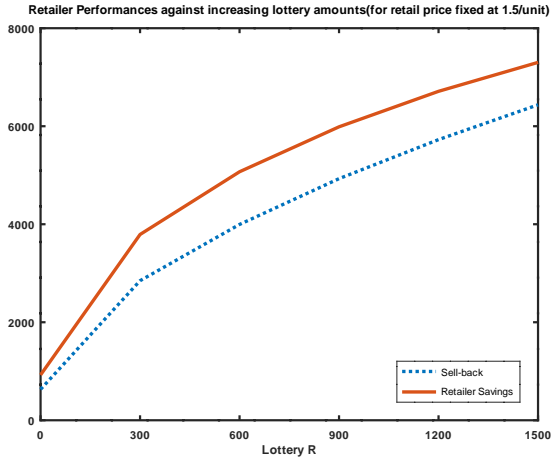}
     \vspace{-0.3in}
     \caption{Variation of metrics against varying lottery prize}
     \label{lottery}
     \vspace{-0.15in}
\end{figure}

We evaluate the performance of a fixed lottery scheme with respect to its incentivization potential and impact on retailer savings (see Fig. \ref{lottery}). Although the probability of a prosumer winning the lottery is very small, humans have a tendency of overweighting small probabilities.  This leads to high sell-back amounts using a reasonable-sized lottery. The retailer savings are also significantly better than in the day-ahead contract pricing scenario, therefore the retailer has a motive to implement the lottery scheme. 

\subsubsection{Effect of Prosumer Penetration}
\begin{figure}[!ht]
     \centering
     \includegraphics[width = 0.48\textwidth]{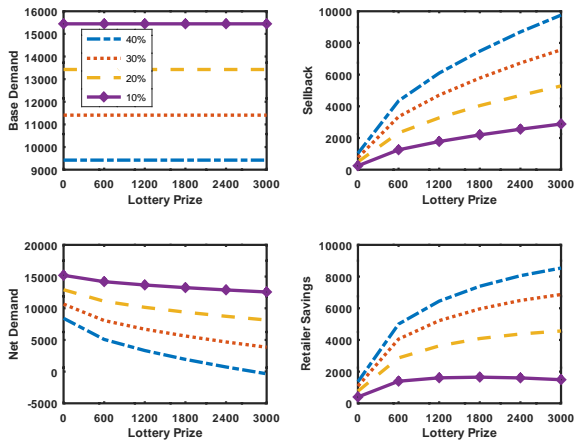}
     \vspace{-0.3in}
     \caption{Variation across metrics for different levels of prosumer penetration}
     \label{penetration}
     \vspace{-0.15in}
\end{figure}

We investigate the effect of prosumer penetration into the community of consumers (see Fig. \ref{penetration}). Specifically, we increase $n$, the number of prosumers from $1000$ to $4000$ while decreasing the number of customers, $N$, from $9,000$ to $6000$. As the prosumer penetration ($\frac{100n \%}{N+n}$) increases, the base demand decreases. For the same lottery prize, the amount of sell-back increases gradually which results in higher savings for the retailer. However, too high penetration is also detrimental to the retailer because the net demand becomes negative, i.e., she has excess energy leftover from sell-back. This would lead to huge loss of revenue for the retailer. 
\section{Conclusions \& Future Work}
In this paper, we first investigated the day-ahead contract pricing for selling back energy to the grid in the lens of prospect theory. We established that the actual sell-back amount can be smaller compared to the amount computed by EUT. Hence, while designing optimal prices, the retailer should consider the CPT approach. We also show that a lottery based incentive scheme may increase the sell-back and savings of the retailer.

In the future, we characterize how the retailer can optimally design the lottery. A prosumer's utility may have temporal correlation. Analyzing and characterization of optimal lottery based scheme under such a scenario also constitutes an interesting future direction. 



\bibliographystyle{IEEEtran}
\bibliography{mybib.bib}

\clearpage
\section{Appendix}
\subsection{Theorem 1}
\begin{proof}
\begin{equation*}
    \begin{multlined}
    E_p(C_k) = p_s C_k + p_e \int_{z_{2,min}}^{C_k}v(z_k-C_k)\frac{d}{dz_k}(\pi(F_z(z_k)))dz_k\\
    = p_s C_k + p_e\int_{z_{2,min}}^{C_k} \lambda (z_k-C_k) \frac{d}{dz_k}(\pi(F_z(z_k)))dz_k \\
    = p_s C_k + \lambda p_e \int_{z_{2,min}}^{C_k}(z_k-C_k) \frac{d}{dz_k}(\pi(F_z(z_k)))dz_k\\
    = p_s C_k + \lambda p_e \int_{z_{2,min}}^{C_k}(z_k-C_k)\pi^{'}(F_z(z_k))f_z(z_k)dz_k
    \end{multlined}
\end{equation*}
\begin{equation*}
    \begin{multlined}
    \frac{dE_p(C_k)}{dC_k} = p_s + \lambda p_e \{(C_k-C_k).\pi^{'}(F_z(C_k))f_z(C_k)\frac{dC_k}{dC_k} \\- (z_{2,min}-C_k)\pi^{'}(F_z(z_{2,min}))f_z(z_{2,min})\frac{dz_{2,min}}{dC_k} +\\ \int_{z_{2,min}}^{C_k}\frac{\partial}{\partial C_k}[(z_k-C_k)\pi^{'}(F_z(z_k))f_z(z_k)]dz_k\}\\
    = p_s + \lambda p_e\int_{z_{2,min}}^{C_k}\frac{\partial}{\partial C_k}[(z_k-C_k)\pi^{'}(F_z(z_k))f_z(z_k)]dz_k\\
    = p_s + \lambda p_e \int_{z_{2,min}}^{C_k} -\pi^{'}(F_z(z_k))f_z(z_k)dz_k\\
    = p_s - \lambda p_e \int_{z_{2,min}}^{C_k} \pi^{'}(F_z(z_k))f_z(z_k)dz_k\\
    = p_s - \lambda p_e\cdot \pi(F_z(z_k))|_{z_{2,min}}^{C_k}\\
    = p_s - \lambda p_e (\pi(F_z(C_k)) - \pi(F_z(z_{2,min}))) = 0 
    \end{multlined}
\end{equation*}
Now, $z_k = s_k - \frac{(\omega_k-p_e)}{\alpha}$. Hence, we have:
\begin{equation*}
    \begin{multlined}
    F_z(z) = P(z_k \leq z) = P(s_k - \frac{(\omega_k-p_e)}{\alpha} \leq z)\\
    = P(s_k \leq z + \frac{(\omega_k-p_e)}{\alpha}) = F_s(z + \frac{(\omega_k-p_e)}{\alpha})
    \end{multlined}
\end{equation*}
This means that $F_z(z_{2,min}) = F_s(s_{min}) = 0$ and $F_z(C_k) = F_s(C_k + \frac{\omega_k-p_e}{\alpha})$. Therefore, we have:
\begin{equation*}
    \begin{multlined}
    \pi(F_s(C_k + \frac{\omega_k-p_e}{\alpha})) = \frac{p_s}{\lambda p_e}\\
    C_k^{*} = F_s^{-1}(\pi^{-1}(\frac{p_s}{\lambda p_e})) - \frac{(\omega_k-p_e)}{\alpha}
    \end{multlined}
\end{equation*}
It is trivial to verify that $C_k^{*}$ is the unique maximizer because $\frac{d^{2}E_p(C_k)}{dC_k^{2}} = -\lambda p_e \pi^{'}(F_z(C_k))f_z(C_k)$ which is negative because $f_z(C_k) > 0$ and $\pi^{'}(F_z(C_k))>0$ (as $\pi(\cdot)$ is monotonically increasing in (0, 1)).
\end{proof}

\subsection{Corollary 1}
\begin{proof} $F_s(\cdot)$ and $\pi(\cdot)$ are both monotonically increasing continuous functions. Therefore, they are \textit{bijective} and hence \textit{invertible}, i.e., each image has a pre-image and the pre-image is unique. Let $f_1$ and $f_2$ be two points in the co-domain of $\pi(\cdot)$ such that $f_1 > f_2$. Since $\pi(\cdot)$ is invertible, $\exists$ $p_1$, $p_2$ $\in$ $\mathcal{D}$($\pi$) such that $p_1 = \pi^{-1}(f_1)$ and $p_2 = \pi^{-1}(f_2)$. Let, if possible, $p_1 \leq p_2$.
\begin{equation*}
    p_1 \leq p_2 \implies \pi(p_1) \leq \pi(p_2) \implies f_1 \leq f_2
\end{equation*}
Thus we arrive at a contradiction. hence, $f_1 > f_2 \implies \pi^{-1}(f_1) > \pi^{-1}(f_2)$. The same can be proved for $F_s(\cdot)$. Therefore, finding the condition for $\pi^{-1}(\frac{p_s}{\lambda p_e}) < \frac{p_s}{p_e}$ is sufficient. 
\begin{equation*}
    \pi^{-1}(\frac{p_s}{\lambda p_e}) < \frac{p_s}{p_e} \implies \frac{p_s}{\lambda p_e} < \pi(\frac{p_s}{p_e})
\end{equation*}
Taking $y = \frac{p_s}{p_e}$, we need to find conditions for $\frac{y}{\pi(y)} < \lambda$. We know that $\lambda > 1$ and $y < \pi(y)$ when $y < \frac{1}{e}$. Therefore, $y < \frac{1}{e} \implies \frac{y}{\pi(y)} < 1 < \lambda$. For higher values of $\lambda$, the size of the interval over which EUT contracts are larger, increases. The proof goes through for any choice of $\pi(\cdot)$ as long as it is monotonically increasing, concave first and then concave with an inflection point in (0, 1) [at some $\frac{1}{r}$ where $r > 1$]. In that case, the results gets modified to the following form : $p_k < \frac{p_e}{r} \implies F_s^{-1}(\pi^{-1}(\frac{p_s}{\lambda p_e})) < F_s^{-1}(\frac{p_s}{p_e})$. 
\end{proof}
\subsection{Corollary 2}
\begin{proof}
We know that $C_{EUT} \geq z_{2,min}$ and $C_{PT} \geq z_{2,min}$. Also, $z_{2,min} \leq z_2 \leq z_{2,max}$. Now, there can be three cases:
\begin{itemize}
\item $z_{2,min} < C_{PT} \leq C_{EUT} < z_1$ : Since both $C_{EUT}$ and $C_{PT}$ are less than $z_1$, according to observation \ref{realtime}, $z^{*}_{EUT} = z^{*}_{PT} = z_1$.

\item $z_{2,min} < C_{PT} < z_1 < C_{EUT}$ :  In this scenario, we have $z^{*}_{PT} = z_1$ and $z^{*}_{EUT} = min(z_2, C_{EUT})$. Since $C_{EUT} > z_1$ and $z_2 > z_1$, $min(z_2, C_{EUT}) > z_1$ which implies $z^{*}_{EUT} > z^{*}_{PT}$.

\item $z_{2,min} < z_1 < C_{PT} \leq C_{EUT}$ : In this scenario, $z_{EUT}^* = \min{(z_2, C_{EUT})} \geq \min{(z_2, C_{PT})} = z_{PT}^*$.
\end{itemize}
Combining all three cases, we conclude that $C_{EUT} \geq C_{PT} \implies z^{*}_{EUT} \geq z^{*}_{PT}$.
\end{proof}



\subsection{Theorem 2}
\begin{proof}
\begin{equation*}
    \begin{multlined}
    U(z) = \omega(s-z)-\frac{\alpha}{2}(s-z)^{2} + R\pi(mz)\\
    = \omega(s-z)-\frac{\alpha}{2}(s-z)^{2} + Rmz + Rg(z)
    = U_1(z) + Rg(z)
    \end{multlined}
\end{equation*}
$g(z) = \pi(mz)-mz$. Since $mz$ is close to $0$, $g(z) > 0$ and $g^{'}(z) > 0$. Also, close to zero, $\pi(\cdot)$ is concave. $g^{''}(z) = m^{2}\pi^{''}(mz) < 0$ which means that $g(\cdot)$ is also concave. $U_1(z)$ is a quadratic function in $z$ and hence has a unique maxima, say at $z^{*}$. Let $z_1$ be the optimal sell-back that maximizes $U(z)$. It can be easily shown that $z_1 > z^{*}$. Therefore, we must have $U_1^{'}(z_1) + Rg^{'}(z_1) = 0$.

Now, let's increase the lottery size infinitesimally to $R^{'}$. Let the optimal sell-back amount now be $z_2$. Hence, $U_1^{'}(z_2) + R^{'}g^{'}(z_2) = 0$.
Let, if possible, $z_2 \leq z_1$. Since the lottery only changes infinitesimally, $z_2 > z^{*}$. Since $U_1(\cdot)$ is quadratic in nature, $U_1^{'}(z_1) < 0$ and $U_1^{'}(z_2) < 0$. Also, $U_1^{'}(z_2) \geq U_1^{'}(z_1)$ since $z^{*} < z_2 \leq z_1$. We already know that $g^{'}(z_2) \geq g^{'}(z_1) > 0$. This means that $U_1^{'}(z_2) + R^{'}g^{'}(z_2) > U_1^{'}(z_1) + Rg^{'}(z_1) > 0$

But, $U_1^{'}(z_2) + R^{'}g^{'}(z_2) = 0$ as $z_2$ is the optimal sell-back for $R^{'}$. Thus, we have a contradiction. Hence, $R^{'} > R \implies z_2 > z_1$. 
\end{proof}
\end{document}